\documentclass[twoside,11pt]{amsart}

\usepackage{amsmath,amssymb,amsbsy,amstext,amsxtra,latexsym}
\usepackage{graphics, epsfig}
\usepackage{hyperref}

\newtheorem{theorem}{Theorem}[section]
\newtheorem{lemma}{Lemma}[section]
\newtheorem{proposition}{Proposition}[section]

\def\mF{{\mathbb F}}
\def\mP{{\mathbb P}}
\def\mQ{{\mathbb Q}}
\def\mC{{\mathbb C}}
\def\tX{{\tilde X}}

\def\cO{{\mathcal O}}
\def\cX{{\mathcal X}}
\def\cY{{\mathcal Y}}

\def\cV{{\mathcal V}}
\def\n{\noindent}
\def\b{\bigskip}
\def\m{\medskip}
\def\t{\tilde}
\def\tX{\tilde X}
\def\T{{\mathbb T}}
\def\Ext{\operatorname{Ext}}
\def\ker{\operatorname{ker}}

\def\Pic{\operatorname{Pic}}
\def\log{\operatorname{log}}
\def\deg{\operatorname{deg}}

\def\Si{{\mathbb F}}

\begin{document}

\title{A Construction of Horikawa Surface via \\
       $\mQ$-Gorenstein Smoothings}

\author{Yongnam Lee and Jongil Park}

\address{Department of Mathematics, Sogang University,
         Sinsu-dong, Mapo-gu, Seoul 121-742, Korea}

\email{ynlee@sogang.ac.kr}

\address{Department of Mathematical Sciences, Seoul National University,
         San 56-1, Sillim-dong, Gwanak-gu, Seoul 151-747, Korea}

\email{jipark@math.snu.ac.kr}

\date{June 1, 2008}



\begin{abstract}
 In this article we prove that Fintushel-Stern's construction of
 Horikawa surface, which is obtained from an elliptic surface
 via a rational blow-down surgery in smooth category,
 can be performed in complex category.
 The main technique involved is $\mQ$-Gorenstein smoothings.
\end{abstract}

\maketitle

\section{Introduction}
\label{sec-1}

\markboth{YONGNAM LEE AND JONGIL PARK}{A CONSTRUCTION OF HORIKAWA
          SURFACE VIA $\mQ$-GORENSTEIN SMOOTHINGS}

 As an application of a rational blow-down surgery on $4$-manifolds,
 R.~Fintushel and R.~Stern showed that Horikawa surface $H(n)$ can
 be obtained from an elliptic surface $E(n)$ via a rational blow-down
 surgery in smooth category~\cite{FS}.
 Note that Horikawa surface $H(n)$ is defined as a double cover of a
 Hirzebruch surface $\mF_{n-3}$ branched over $|6C_0+(4n-8)f|$,
 where $C_0$ is a negative section and $f$ is a fiber of $\mF_{n-3}$.

 In this article we show that a rational blow-down surgery to obtain
 Horikawa surface can be performed in fact in complex category.
 That is, we reinterpret algebraically Fintushel-Stern's topological
 construction~\cite{FS} of Horikawa surface $H(n)$ to give a complex
 structure on it.
 The main technique we use in this paper is $\mQ$-Gorenstein smoothings.
 Note that $\mQ$-Gorenstein smoothing theory developed in deformation
 theory in last thirty years is a very powerful tool
 to construct a non-singular surface of general type.
 The basic scheme is the following: Suppose that a projective surface
 contains several disjoint chains of curves representing the resolution
 graphs of special quotient singularities.
 Then, by contracting these chains of curves,
 we get a singular surface $X$ with special quotient singularities.
 And then we investigate the existence of a $\mQ$-Gorenstein
 smoothing of $X$. It is known that the cohomology $H^2(T^0_X)$
 contains the obstruction space of a $\mQ$-Gorenstein smoothing of $X$.
 That is,
 if $H^2(T^0_X)=0$, then there is a $\mQ$-Gorenstein smoothing of $X$.
 For example,
 we recently constructed a simply connected minimal surface of general type
 with $p_g=0$ and $K^2 =2$ by proving the cohomology $H^2(T^0_X)=0$~\cite{LP}.
 But, in general, the cohomology $H^2(T^0_X)$ is not zero
 and it is a very difficult problem to determine
 whether there exists a $\mQ$-Gorenstein smoothing of $X$.
 In this article we also give a family of examples which admit
 $\mQ$-Gorenstein smoothings
 even though the cohomology $H^2(T^0_X)$ does not vanish.
 Our main technique is a $\mQ$-Gorenstein smoothing theory with
 a cyclic group action. It is briefly reviewed and developed in Section~2.

 The sketch of our construction whose details are given in Section~3
 is as follows:
 We first construct a simply connected relatively minimal elliptic surface
 $E(n)$ ($n \ge 5$) with a special fiber,
 which contains two linear chains of configurations of $\mP^1$'s
 \[\underset{U_{n-3}}{\overset{-n}{\circ}}-\underset{U_{n-4}}
   {\overset{-2}{\circ}} -\underset{U_{n-5}}{\overset{-2}{\circ}}-
   \cdots - \underset{U_1}{\overset{-2}{\circ}}.\]

\n
 We construct this kind of an elliptic surface $E(n)$ explicitly
 by using a double cover of a blowing-up of Hirzebruch surface $\mF_{n}$
 branched over a special curve. The double cover of $\mF_{n}$ has
 two rational double points $A_1$ and $A_{2n-9}$.
 Then its minimal resolution is an elliptic surface $E(n)$
 which has an $I_{2n-6}$ as a special fiber. Now we contract
 these two linear chains of configurations of $\mP^1$'s to produce a
 normal projective surface $X_n$ with two special quotient
 singularities, both singularities are of type $\frac{1}{(n-2)^2}(1, n-3)$.
 Finally, we apply $\mQ$-Gorenstein smoothing theory with
 a cyclic group action developed in Section~2 for $X_n$
 in order to get our main result which is following.

\begin{theorem}
\label{thm-main}
 The projective surface $X_n$ obtained by contracting two disjoint
 configurations $C_{n-2}$ from an elliptic surface $E(n)$ admits
 a $\mQ$-Gorenstein smoothing of two quotient singularities all together,
 and a general fiber of the $\mQ$-Gorenstein smoothing
 is Horikawa surface $H(n)$.
\end{theorem}

\m

\n{\em Acknowledgements}.
 The authors would like to thank Ronald Fintushel and Ronald Stern
 for helpful comments on Horikawa surfaces during the Conference of
 Algebraic Surfaces and $4$-manifolds held at KIAS. The authors also
 wish to thank Roberto Pignatelli for valuable discussions to prove
 Proposition 3.1. Yongnam Lee was supported by Korea Research Foundation
 Grant funded by the Korean Government (KRF-2005-070-C00005).
 He would like to thank Meng Chen for his generous hospitality during
 his visit Fudan University, where some of the paper was completed.
 Jongil Park was supported by SBS Foundation Grant in 2007 and
 he also holds a joint appointment in the Research
 Institute of Mathematics, Seoul National University.

\b

\section{$\mQ$-Gorenstein smoothing}
\label{sec-2}

 In this section we briefly review a theory of $\mQ$-Gorenstein smoothing
 for projective surfaces with special quotient singularities,
 which is a key technical ingredient in our main construction.

\m

 \n {\bf Definition.} Let $X$ be a normal projective surface with quotient
 singularities. Let $\cX\to\Delta$ (or $\cX/\Delta$) be a flat family of
 projective surfaces over a small disk $\Delta$. The one-parameter
 family of surfaces $\cX\to\Delta$ is called a {\it
 $\mQ$-Gorenstein smoothing} of $X$ if it satisfies the following
 three conditions;

\n(i) the general fiber $X_t$ is a smooth projective surface,

\n(ii) the central fiber $X_0$ is $X$,

\n(iii) the canonical divisor $K_{\cX/\Delta}$ is $\mQ$-Cartier.

\m

 A $\mQ$-Gorenstein smoothing for a germ of a quotient singularity
 $(X_0, 0)$ is defined similarly. A quotient singularity which
 admits a $\mQ$-Gorenstein smoothing is called a {\it singularity
 of class T}.

\begin{proposition}[\cite{KSB, Man91, Wa1}]
\label{pro-2.1}
 Let $(X_0, 0)$ be a germ of two dimensional quotient
 singularity. If $(X_0, 0)$ admits a $\mQ$-Gorenstein smoothing over
 the disk, then $(X_0, 0)$ is either a rational double point or a
 cyclic quotient singularity of type $\frac{1}{dn^2}(1, dna-1)$ for some
 integers $a, n, d$ with $a$ and $n$ relatively prime.
\end{proposition}

\begin{proposition}[\cite{KSB, Man91, Wa2}]
\label{pro-2.2}
\begin{enumerate}
 \item The singularities ${\overset{-4}{\circ}}$ and
 ${\overset{-3}{\circ}}-{\overset{-2}{\circ}}-{\overset{-2}{\circ}}-\cdots-
 {\overset{-2}{\circ}}-{\overset{-3}{\circ}}$ are of class $T$.
 \item If the singularity
 ${\overset{-b_1}{\circ}}-\cdots-{\overset{-b_r}{\circ}}$ is of
 class $T$, then so are
 $${\overset{-2}{\circ}}-{\overset{-b_1}{\circ}}-\cdots-{\overset{-b_{r-1}}
 {\circ}}- {\overset{-b_r-1}{\circ}} \quad\text{and}\quad
 {\overset{-b_1-1}{\circ}}-{\overset{-b_2}{\circ}}-\cdots-
 {\overset{-b_r}{\circ}}-{\overset{-2}{\circ}}.$$ \item Every
 singularity of class $T$ that is not a rational double point can be
 obtained by starting with one of the singularities described in
 $(1)$ and iterating the steps described in $(2)$.
\end{enumerate}
\end{proposition}

 Let $X$ be a normal projective surface with singularities of class
 $T$. Due to the result of Koll\'ar and Shepherd-Barron \cite{KSB},
 there is a $\mQ$-Gorenstein smoothing locally for each singularity
 of class $T$ on $X$ (see Proposition~\ref{pro-ksb}).
 The natural question arises whether this local $\mQ$-Gorenstein
 smoothing can be extended over the global surface $X$ or not.
 Roughly geometric interpretation is the following:
 Let $\cup_{\alpha} V_{\alpha}$ be an open covering of $X$ such that
 each $V_{\alpha}$ has at most one singularity of class $T$.
 By the existence of a local $\mQ$-Gorenstein
 smoothing, there is a $\mQ$-Gorenstein smoothing
 $\cV_{\alpha}/\Delta$. The question is if these families glue to
 a global one. The answer can be obtained by figuring out the
 obstruction map of the sheaves of deformation
 $T^i_X=Ext^i_X(\Omega_X,\cO_X)$ for $i=0,1,2$.
 For example, if $X$ is a smooth surface,
 then $T^0_X$ is the usual holomorphic tangent sheaf $T_X$ and
 $T^1_X=T^2_X=0$. By applying the standard result of deformations
 \cite{LS, Pal} to a normal projective surface with quotient
 singularities, we get the following

\begin{proposition}[\cite{Wa1}, \S 4]
\label{pro-2.3}
 Let $X$ be a normal
 projective surface with quotient singularities. Then
\begin{enumerate}
 \item The first order deformation space of $X$ is represented by
 the global Ext 1-group $\T^1_X=\Ext^1_X(\Omega_X, \cO_X)$. \item
 The obstruction lies in the global Ext 2-group
 $\T^2_X=\Ext^2_X(\Omega_X, \cO_X)$.
 \end{enumerate}
\end{proposition}

 Furthermore, by applying the general result of local-global spectral
 sequence of ext sheaves (\cite{Pal}, \S 3) to deformation theory
 of surfaces with quotient singularities so that
 $E_2^{p, q}=H^p(T^q_X) \Rightarrow \T^{p+q}_X$,
 and by $H^j(T^i_X)=0$ for $i, j\ge 1$, we also get

\begin{proposition}[\cite{Man91, Wa1}]
\label{pro-2.4}
 Let $X$ be a normal projective surface with quotient singularities.
 Then
 \begin{enumerate}
 \item We have the exact sequence
 $$0\to H^1(T^0_X)\to \T^1_X\to \ker [H^0(T^1_X)\to H^2(T^0_X)]\to 0$$
 where $H^1(T^0_X)$ represents the first order deformations of $X$
 for which the singularities remain locally a product.
 \item If $H^2(T^0_X)=0$, every local deformation of
 the singularities may be globalized.
 \end{enumerate}
\end{proposition}

 The vanishing $H^2(T^0_X)=0$ can be obtained via
 the vanishing of $H^2(T_V(-\log \ E))$, where $V$ is the minimal resolution
 of $X$ and $E$ is the reduced exceptional divisors. Note that
 every singularity of class $T$ has a local $\mQ$-Gorenstein
 smoothing by Proposition~\ref{pro-ksb} below.

 Let $X$ be a normal projective surface with singularities of class $T$.
 Our concern is to understand $\mQ$-Gorenstein smoothings in $\T^1_X$,
 not the whole first order deformations. These special deformations
 can be constructed via local index one cover. Let $U\subset X$ be an
 analytic neighborhood with an index one cover $U'$. For the case of
 the field $\mC$, this index one cover is unique up to isomorphism.
 The first order deformations which associate $\mQ$-Gorenstein
 smoothings can be realized as the invariant part of $T^1_{U'}$.
 The sheaves $\t T^1_X$ are defined by the index one covering stack
 and by the \'etale sites \cite{Hac}. The first order deformation of
 a $\mQ$-Gorenstein smoothing of singularities of class $T$ is
 expressed by the cohomology $H^0(\tilde T^1_X)$ \cite{Hac, Has, KSB}.
 By the help of the birational geometry in threefolds and their
 applications to deformations of surface singularities,
 the following proposition is obtained.
 Note that the cohomology $H^0(\t{T}^1_X)$ is given explicitly as follows.

\begin{proposition}[\cite{KSB, Man91}]
\label{pro-ksb}
\begin{enumerate}
 \item Let $a, d, n>0$ be integers with $a, n$ relatively prime and
 consider a map $\pi :\cY/\mu_n \to \mC^d$, where $\cY\subset
 \mC^3\times \mC^d$ is the hypersurface of equation
 $uv-y^{dn}=\sum_{k=0}^{d-1}t_ky^{kn};\ t_0,\ldots, t_{d-1}$ are
 linear coordinates over $\mC^d$, $\mu_n$ acts on $\cY$ by
 $$\mu_n\ni\xi : (u, v, y, t_0,\ldots, t_{d-1})\to (\xi u, \xi^{-1}v,
 \xi^ay, t_0, \ldots, t_{d-1})$$ and $\pi$ is the factorization to
 the quotient of the projection $\cY\to\mC^d$. Then $\pi$ is a
 $\mQ$-Gorenstein smoothing of the cyclic singularity of a germ
 $(X_0, 0)$ of type $\frac{1}{dn^2}(1, dna-1)$. Moreover every
 $\mQ$-Gorenstein smoothing of $(X_0, 0)$ is isomorphic to the
 pull-back of $\pi$ for some germ of holomorphic map
 $(\mC, 0)\to (\mC^d, 0)$.
 \item Let $X$ be a normal projective surface with singularities of class T.
 Then
 $$H^0(\t{T}^1_X)=\sum_{p\,\in\text{\,singular points of $X$}}
 \mC_p^{\oplus d_p}$$
 where a singular point $p$ is of type $\frac{1}{d_p n^2}
 (1, d_p an-1)$ with $(a, n)=1$.
\end{enumerate}
\end{proposition}

\begin{theorem}[\cite{LP}]
\label{thm-2.1}
 Let $X$ be a normal projective surface with
 singularities of class $T$.
 Let $\pi: V\to X$ be the minimal resolution and
 let $E$ be the reduced exceptional
 divisors. Suppose that $H^2(T_V(-\log \ E))=0$.
 Then $H^2(T^0_X)=0$ and there is a $\mQ$-Gorenstein smoothing of $X$.
\end{theorem}

 As we see in Theorem~\ref{thm-2.1} above,
 if $H^2(T^0_X)=0$, then there is a $\mQ$-Gorenstein smoothing of $X$.
 For example,
 we constructed a simply connected minimal surface of general type with
 $p_g=0$ and $K^2 =2$ by proving the cohomology $H^2(T^0_X)=0$~\cite{LP}.
 But, in general, the cohomology $H^2(T^0_X)$ is not zero
 and it is a very difficult problem to determine
 whether there exists a $\mQ$-Gorenstein smoothing of $X$.
 Hence, in the case that $H^2(T^0_X) \not =0$,
 we have to develop another technique in order to investigate the existence
 of  $\mQ$-Gorenstein smoothings.
 Even though we do not know whether such a technique exists in general,
 if $X$ is a normal projective surface with
 singularities of class $T$ which admits a cyclic group
 with some nice properties, then we are able to show that it admits
 a $\mQ$-Gorenstein smoothing. Explicitly, we get the following theorem.

\begin{theorem}
\label{thm-2.2}
 Let $X$ be a normal projective surface with
 singularities of class $T$. Assume that a cyclic group
 $G$ acts on $X$ such that
 \begin{enumerate}
 \item $Y=X/G$ is a normal projective
 surface with singularities of $T$,
 \item  $p_g(Y)=q(Y)=0$,
 \item $Y$ has a $\mQ$-Gorenstein smoothing,
 \item the map $\sigma: X\to Y$ induced by a cyclic covering is flat,
  and the branch locus $D$ (resp. the ramification locus) of the map
  $\sigma: X\ \to Y$ is an irreducible nonsingular curve lying outside
  the singular locus of $Y$ (resp. of $X$), and
 \item $H^1(Y, \cO_Y(D))=0$.
\end{enumerate}
 Then there exists a $\mQ$-Gorenstein smoothing of $X$ that is compatible
 with a $\mQ$-Gorenstein smoothing of $Y$.
 And the cyclic covering extends to the $\mQ$-Gorenstein smoothing.
\end{theorem}

\begin{proof}
 Let $\cY\to\Delta$ be a $\mQ$-Gorenstein smoothing of $Y$, and let
 $Y_t$ be a general fiber of the $\mQ$-Gorenstein smoothing.
 By the semi-continuity, we have $p_g(Y_t)=q(Y_t)=0$.
 The base change theorem and Leray spectral
 sequence imply that $H^1(\cY, \cO_{\cY})=H^2(\cY, \cO_{\cY})=0$.
 It gives an isomorphism $r_0: \Pic (\cY)\to \Pic(Y)$ and an
 injective map $r_t: \Pic(\cY)\to \Pic(Y_t)$ (Lemma 2 in~\cite{Man91}).
 The vanishing $H^1(Y, \cO_Y(D))=0$ ensures that the deformation of
 $Y$ can be lifted to the deformation of the pair $(Y, D)$, i.e.
 the branch divisor $D$ is extended to $D_t$ in $Y_t$. Since the divisor
 $D$ is nonsingular, $D_t$ is also nonsingular.
 And the flatness of the map ensures that the divisor $L$
 which is the data of the cyclic cover, i.e. $L^{\otimes |G|}\cong D$,
 is extended to $L_t$ with $L_t^{\otimes |G|}\cong D_t$.
 Hence, the cyclic covering extends to the $\mQ$-Gorenstein smoothing of $Y$.
\end{proof}

\b

\section{A construction of Horikawa surface}
\label{sec-3}

 Let $E(n)$ be a simply connected relatively minimal
 elliptic surface with a section and with $c_2=12n$.
 Then
 there is only one up to diffeomorphism such an elliptic surface and
 the canonical class is given by
 $K_{E(n)}=(n-2)C$, where $C$ is a general fiber of an elliptic fibration.
 Hence each section is a nonsingular rational curve whose
 self-intersection number is $-n$.
 Assume that $n \ge 4$ and
 let $C_{n-2}$ denote a simply connected smooth $4$-manifold
 obtained by pluming the $(n-3)$ disk bundles over the $2$-sphere
 according to the linear diagram
\[\underset{U_{n-3}}{\overset{-n}{\circ}}-\underset{U_{n-4}}
   {\overset{-2}{\circ}} -\underset{U_{n-5}}{\overset{-2}{\circ}}-
   \cdots - \underset{U_1}{\overset{-2}{\circ}}.\]

 Assume that an elliptic surface $E(n)$ has two configurations $C_{n-2}$
 such that all embedded $2$-spheres $U_i$ are holomorphic curves
 (We show the existence of such an $E(n)$ later).
 Let $Y_n'$ be a normal projective surface obtained by contracting
 one configuration $C_{n-2}$ from $E(n)$. Then $Y_n'$ does not admit a
 $\mQ$-Gorenstein smoothing  because it violates Noether inequality
 (Corollary 7.5 in~\cite{FS}).
 In fact, it does not satisfy the vanishing condition in the hypothesis of
 Theorem~\ref{thm-2.1}, that is, we have $H^2(E(n), T_{E(n)})\ne 0$:
 Let $h: E(n)\to \mP^1$ be an elliptic fibration.
 Assume that $C$ is a general fiber of the map $h$.
 We have an injective map $0\to h^*\Omega_{\mP^1} \to \Omega_{E(n)}$
 and the map induces an injection
 $H^0(\mP^1,\Omega_{\mP^1}(n-2))\hookrightarrow H^0(E(n),\Omega_{E(n)}((n-2)C))$
 by tensoring $(n-2)C$ on $0\to h^*\Omega_{\mP^1} \to \Omega_{E(n)}$.
 Since $K_{E(n)}=(n-2)C$,
 the cohomology $H^0(E(n), \Omega_{E(n)}(K_{E(n)}))$ is not zero.
 Hence the Serre duality implies that $H^2(E(n), T_{E(n)})$ is not zero.

 Next, let $X_n$ be a normal projective surface obtained by contracting
 two disjoint configurations $C_{n-2}$ from $E(n)$, and we want
 to investigate the existence of a $\mQ$-Gorenstein smoothing of $X_n$.
 As a warming-up, we first investigate $n=4$ case.

\m

\n {\bf Example}. R. Gompf constructed a family of symplectic
 $4$-manifolds by taking a fiber sum of other symplectic
 $4$-manifolds~\cite{Gom}.
 To recall Gompf's example briefly, we start with a simply connected
 relatively minimal elliptic surface $E(4)$ with a section and with $c_2=48$.
 It is known that $E(4)$ admits nine rational $(-4)$-curves as disjoint sections.
 Rationally blowing-down $n$ $(-4)$-curves of $E(4)$ is the same as
 the normal connected sum of $E(4)$ with $n$ copies of $\mP^2$
 by identifying a conic in each $\mP^2$ with one $(-4)$-curve in $E(4)$.
 Let us denote this $4$-manifold by $W_{4,n}$.
 Then the manifold $W_{4, 1}$ does not admit any complex structure
 because it violates the Noether inequality.
 But we will show that $W_{4,2}$ admits a complex structure using
 a $\mQ$-Gorenstein smoothing theory.
 For this, let us first denote the singular projective surface obtained by
 contracting $n$ $(-4)$-sections from $E(4)$ by $W_{4,n}'$.
 And then we claim that $W_{4,2}'$ has a $\mQ$-Gorenstein smoothing.
 The reason is following:
 Consider $E(4)$ as a double cover of Hirzebruch surface $\mF_4$ branched over
 an irreducible nonsingular curve $D$ in the linear system $|4(C_0+4f)|$,
 where $C_0$ is a negative section and $f$ is a fiber of $\mF_4$.
 Then $H^1(\mF_4, \cO_{\mF_4}(D))=0$: Since $p_g(\mF_4)=q(\mF_4)=0$,
 \[H^1(\mF_4, \cO_{\mF_4}(D))\simeq H^1(D, \cO_D(D))\simeq
   H^0(D, \cO_D(K_D-D))^{\vee}.\]
 And $\deg K_D-D^2=4(C_0+4f)(2C_0+10f)-16(C_0+4f)^2=-24<0$ implies that
 $H^0(D, \cO_D(K_D-D))=0$.
 Since $D$ does not intersect $C_0$, $W_{4,2}'$ is a double cover
 of a cone $\hat \mF_4$ which is a contraction of $C_0$ from $\mF_4$.
 This implies that the map $\sigma$ induced by a double cover is flat
 and $H^1(Y, \cO_Y(D))=0$.
 Note that $\hat \mF_4$ has a $\mQ$-Gorenstein smoothing whose general fiber
 is $\mP^2$. It is obtained by a pencil of hyperplane section of the cone
 of the Veronese surface imbedded in $\mP^5$.
 Hence $W_{4, 2}'$ has a $\mQ$-Gorenstein smoothing by Theorem~\ref{thm-2.2}.
 Finally, since the rational blow-down manifold $W_{4, 2}$ is
 diffeomorphic to the general fiber of the $\mQ$-Gorenstein smoothing
 of $W_{4, 2}'$, $W_{4, 2}$ admits a complex structure.
 Furthermore, using a triple
 cover of $\Si_4$ branched over $D$ in the linear system $|3(C_0+4f)|$,
 we can also prove that $W_{4, 3}'$ has a $\mQ$-Gorenstein smoothing
 by the similar proof as above.
 And, by extending Theorem~\ref{thm-2.2} to a finite abelian group,
 it is possible to show that some other manifolds $W_{4,n}'$
 has a $\mQ$-Gorenstein smoothing, too.
 We leave it for a future research.

 \m

 Now we investigate the general case. Assume that $n \ge 5$ and
 let $\mF_n$ be a Hirzebruch surface. Let $C_0$ be
 a negative section with $C_0^2=-n$ and $f$ be a fiber of $\mF_n$.
 Consider the linear system $|4(C_0+nf)|$.
 The surface $\mF_n$ can be obtained from the cone over
 a rational normal curve of degree $n$ by blowing up the vertex.
 And a curve in the linear system $|4(C_0+nf)|$ is the strict transform
 of the hyperplane section of the cone.
 By Bertini's theorem, there is an irreducible
 nonsingular curve in the linear system $|4(C_0+nf)|$.
 The double cover of $\mF_n$ branched over an irreducible
 nonsingular member in $|4(C_0+nf)|$ is an elliptic surface $E(n)$:
 Let $\sigma: \hat X_n\to\mF_n$ be a double covering branched over
 an irreducible nonsingular member in the linear system $|4(C_0+nf)|$.
 Then, by the invariants of a double covering (\cite{BHPV}, Chapter V),
 we have
 $p_g(\hat X_n)=p_g(\mF_n)+h^0(\mF_n, K_{\mF_n}+L)=h^0(\mF_n, (n-2)f)=n-1$
 and
 $\chi(\cO_{\hat X_n})=2\chi(\cO_{\mF_n})+\frac{1}{2}(L\cdot
 K_{\mF_n})+\frac{1}{2}(L\cdot L)=n$, where $L=2(C_0+nf)$.
 Therefore we have $q(\hat X_n)=0$ and
 $K_{\hat X_n}^2=2(\sigma^*(K_{\mF_n}+L))^2=2((n-2)f)^2=0$.

\m

 In this article we want to choose a special irreducible (singular) curve
 $D$ in the linear system $|4(C_0+nf)|$, which has a special intersection
 with one special fiber $f$: Note that $D\cdot f=4$.
 We want $D$ to intersect with $f$ at two distinct points $p$ and
 $q$ that are not in $C_0$.
 Let $x=0$ be the local equation of $f$ and $x, y$ be a coordinate
 at $p$ (resp. at $q$). We require that the local equation of $D$ at $p$
 (resp. at $q$) is $(y-x)(y+x)=0$ (resp. $(y-x^{n-4})(y+x^{n-4})=0$).
 These are $3(n-4)+3$-conditions: $1, x, x^2, \ldots, x^{2n-9}, y,
 yx, \ldots, yx^{n-5}$ terms should vanish to have the local
 analytic equation $(y-x^{n-4})(y+x^{n-4})=0$.
 By next lemmas and proposition,
 we have such a curve $D$ satisfying the conditions above.

 \begin{lemma}
\label{lem-3.1} We have $h^0(\mF_n, \cO_{\mF_n}(D))=10n+5$, where
$D$ is a member in the linear system $|4(C_0+nf)|$.
\end{lemma}

 \begin{proof}
 Let $C=C_0+nf$. Then, by the following two exact sequences
 \[ 0\to\cO_{\mF_n}\to \cO_{\mF_n}(f)\to \cO_{\mP^1}\to 0\]
 \[ 0\to\cO_{\mF_n}(kf)\to \cO_{\mF_n}((k+1)f)\to \cO_{\mP^1}\to 0, \]
 we have $h^1(\mF_n, \cO_{\mF_n}(kf))=0$ for all nonnegative integers $k$.
 And from the exact sequence
 \[ 0\to\cO_{\mF_n}(nf)\to \cO_{\mF_n}(C)\to \cO_{C_0}(C)\to 0,\]
 we also have
 $h^1(\mF_n, \cO_{\mF_n}(C))=0$ and $h^0(\mF_n, \cO_{\mF_n}(C))=n+2$.
 Hence, by considering the exact sequences similarly
 \[0\to \cO_{\mF_n}(kC)\to\cO_{\mF_n}((k+1)C)\to
 \cO_C((k+1)C)\to 0 ,\]
 we finally get $h^0(\mF_n, \cO_{\mF_n}(D))=10n+5$.
 \end{proof}

\m

 First, we assume that $D$ is nonsingular at every point except the two points
 $p$ and $q$.
 Let $\sigma: \tX_n\to\mF_n$ be a double covering branched over
 the curve $D$ chosen above. Then $\tX_n$ is a singular elliptic surface
 with $p_g=n-1$ and $\chi =n$
 which has two rational double points by the local equations of $D$
 at $p$ and $q$ - one is $A_1$ ($z^2=y^2-x^2$) and the other one is
 $A_{2n-9}$ ($z^2=y^2-x^{2n-8}$).
 Therefore its minimal resolution is also an elliptic surface $E(n)$.
 First we blow up at $p$ and $q$ in $\mF_n$. Then we have an
 exceptional curve coming from a blowing up at $p$ which intersects
 with the proper transform of $D$ transversally at two points,
 and we also have an exceptional curve coming from a blowing up at $q$
 which intersects with the proper transform of $D$ at one point, say $q_1$.
 Let $x=0$ be the local equation of the $(-1)$-exceptional curve at
 $q_1$. Then the local equation of the proper transform of $D$ at
 $q_1$ is $(y-x^{n-5})(y+x^{n-5})=0$. We blow up again at $q_1$. By
 the continuation of blowing up at infinitely near points of $q$,
 we have the following configuration of smooth rational curves

 \[\begin{array}{cccc}
   \underset{U_{n-3}}{\overset{-n}{\circ}} - &
   \underset{U_{n-4}}{\overset{-2}{\circ}} &
   - \underset{U_{n-5}}{\overset{-2}{\circ}}-
   \cdots - &
  \underset{U_1}{\overset{-2}{\circ}}\\
   & \vert & & \vert \\
   & {\underset{E_1, -1}{\circ}} & & {\underset{E_2, -1}{\circ}}
  \end{array}\]
 where the proper transform of $D$ intersects with $E_i$, $i=1,2$
 at two points transversally. We denote this surface by $Z_n$,
 which is obtained by $(n-3)$ times blowing-ups of $\mF_n$.

 \m

 Let $\pi: Z_n\to\mF_n$ be a map and
 $\Delta=\pi^*(4C_0+4nf)-2E_1-2U_{n-5}-4U_{n-6}-6U_{n-7}-
 \cdots-2(n-5)U_1-2(n-4)E_2$.
 For a simple computation, we write it as
 \[\Delta=\pi^*(4C_0+4nf)-2F-2\sum_{i=1}^{n-4}F_i ,\]
 where $F=E_1,\, F_1=U_{n-5}+U_{n-6}+\cdots +E_2,\, F_2=U_{n-6}+\cdots +E_2,
 \ldots, F_{n-4}=E_2$.
 Note that $F_i$ is not necessarily irreducible and
 $F^2=F_i^2=K_{Z_n}\cdot F=K_{Z_n}\cdot F_i=-1$ for all $i=1, \ldots, n-4$.
 Let $f_0=U_{n-4}$, which is a proper transform of the fiber,
 and let $L=\Delta-(\pi^*C_0+f_0)-K_{Z_n}$.

\m

 In Proposition 3.1 below,
 we prove that the linear system $|\Delta|$ is base point free.
 Then, by Bertini's theorem, we conclude that
 $D$ is nonsingular except the two points $p$ and $q$.

\begin{lemma}
\label{lem-3.2}
 $L^2\ge 5$ and $L$ is nef.
\end{lemma}

\begin{proof}
 Since $\pi^*f=f_0+F+\sum_{i=1}^{n-4}F_i$, we have
 $L=\pi^*(5C_0+(5n+2)f)-f_0-3F-3\sum_{i=1}^{n-4}F_i=\pi^*(5C_0+(5n-1)f)+2f_0$.
 Furthermore, since the two points $p$ and $q$ are not in $C_0$,
 we also have $C_0\cdot f=C_0\cdot f_0=1$ and $f_0^2=-2$.
 Therefore
 \[L^2=-25n+(5n-1)10+20-8=25n+2\ge 5.\]
 Let $G$ be an irreducible curve which is neither $f_0$ nor $C_0$.
 Note that $L\cdot f_0=5-4=1$ and $L\cdot C_0=-5n+5n-1+2=1$. Write
 $\pi_*G=aC_0+bf$. Then we have
 \[G\cdot L\ge \pi_*G \cdot(5C_0+(5n-1)f)=(aC_0+bf)\cdot(5C_0+(5n-1)f)=-a+5b.\]
 We note that the linear system $|aC_0+bf|$ contains an irreducible
 curve in $\mF_n$ if and only if $a=0, b=1$; or $a=1, b=0$; or $a>0,
 b>an$; or $a>0, b=an$ with $n>0$
 (refer to Corollary 2.18, Chapter V in~\cite{Hart}).
 Therefore it is impossible that $G\cdot L<0$.
 It implies that $L$ is nef.
\end{proof}

\begin{lemma}
\label{lem-3.3}
 The linear system $|\Delta-(\pi^*C_0+f_0)|$ on $Z_n$
 is base point free.
\end{lemma}

\begin{proof}
 By Lemma~\ref{lem-3.2} above, $L$ is nef and $L^2\ge 5$.
 Hence, applying to Reider's theorem \cite{Rdr}, if the adjoint linear
 series $|\Delta-(\pi^*C_0+f_0)|=|K_{Z_n}+L|$ has a base point at $x$
 then there is an effective divisor $G$ in $Z_n$ passing through $x$
 such that either $G\cdot L=0$ and $G^2=-1$; or $G\cdot L=1$ and
 $G^2=0$.

 Assume that $G\cdot L=0$ and $G^2=-1$. Write $G=G_1+\cdots +G_k$,
 where $G_k$ is an irreducible curve. Since $L$ is nef,
 $G\cdot L=0$ implies that $G_i \cdot L=0$ for all $i=1, \ldots, k$.
 Then we get a contradiction by a similar argument to show that
 $-a+5b\le 0$ in the proof Lemma~\ref{lem-3.2} above.

 Assume that $G\cdot L=1$ and $G^2=0$. By the same argument in the case
 $G\cdot L=0$, $G=G_1$. Then we get a contradiction by a similar argument
 to show that $-a+5b\le 1$ unless $G=C_0$ or $f_0$.
 Furthermore, since $C_0^2=-n$ and $f_0^2=-2$, it also contradicts.
\end{proof}

\begin{proposition}
\label{pro-3.1}
 The linear system $|\Delta|$ on $Z_n$ is base point free.
\end{proposition}

\begin{proof}
 Note that $\Delta\cdot \pi^*C_0=\Delta\cdot f_0=0$.
 Therefore we have $\cO_{\pi^*C_0+f_0}(\Delta)=\cO_{\pi^*C_0+f_0}$.
 Furthermore, by Lemma~\ref{lem-3.2} above and by the vanishing theorem,
 we also have
 \[H^1(Z_n, \Delta-(\pi^*C_0+f_0))=H^1(Z_n, K_{Z_n}+L)=0.\]
 Hence, using Lemma~\ref{lem-3.3} above and the short exact sequence
 \[0\to \cO_{Z_n}(\Delta-(\pi^*C_0+f_0))\to\cO_{Z_n}(\Delta)
  \to\cO_{\pi^*C_0+f_0}\to 0 ,\]
 we conclude that the linear system $|\Delta|$ is base point free.
\end{proof}

 Next, by Artin's criterion of contraction~\cite{Art},
 we can contract a configuration $C_{n-2}$,
 which is a linear chain of $\mP^1$'s
 \[\underset{U_{n-3}}{\overset{-n}{\circ}}-\underset{U_{n-4}}
   {\overset{-2}{\circ}} -\underset{U_{n-5}}{\overset{-2}{\circ}}-
   \cdots - \underset{U_1}{\overset{-2}{\circ}} , \]
 so that it produces a singular normal projective surface.
 We denote this surface by $Y_n$.
 We note that $\Delta$ is the proper transform of $D$ in $Z_n$ and
 that $Y_n$ has a cyclic quotient singularity of type
 $\frac{1}{(n-2)^2}(1, n-3)$, which is a singularity of class $T$.

\begin{lemma}
\label{lem-3.4}
 $H^1(Z_n, \cO_{Z_n}(\Delta))=0$.
\end{lemma}

\begin{proof}
 Since $p_g(Z_n)=q(Z_n)=0$, we have
 \[H^1(Z_n, \cO_{Z_n}(\Delta))\simeq H^1(\Delta, \cO_{\Delta}(\Delta))
 \simeq H^0(\Delta, \cO_{\Delta}(K_{\Delta}-{\Delta}))^{\vee}.\]
 We also have ${\Delta}^2=D^2-4(n-3)=12n+12$ and $\deg K_{\Delta}=\deg
 K_D-2(n-3)=12n-8-2(n-3)=10n-2$. Therefore it satisfies
 \[\deg K_{\Delta}-{\Delta}^2=10n-2-12n-12<0 ,\]
 and it implies that $H^0(\Delta, \cO_{\Delta}(K_{\Delta}-\Delta))=0$.
\end{proof}

\begin{proposition}
\label{pro-3.2}
  The singular surface $Y_n$ admits a $\mQ$-Gorenstein smoothing.
\end{proposition}

\begin{proof}
 It is enough to show that
 $-K_{Y_n}$ is effective (Theorem 21 in~\cite{Man91}).
 Let $\pi: Z_n\to\mF_n$ be a composition of blowing-ups,
 and $\psi: Z_n\to Y_n$ be a contraction. Then we have
 \[K_{Z_n}=\pi^*K_{\mF_n}+E_1+U_{n-5}+2U_{n-6}+\cdots +(n-5)U_1+(n-4)E_2.\]
 Since $K_{\mF_n}=-2C_0-(n+2)f$, $-K_{Z_n}$ is effective.
 Furthermore, since $h^0(-K_{Y_n})=h^0(\psi_*(-K_{Z_n}))= h^0(-K_{Z_n})$
 (\S  3.9.2 in~\cite{Wa1}), $-K_{Y_n}$ is also effective.
\end{proof}

\m

 Now we are in a position to prove our main theorem mentioned
 in the Introduction.
 First remind that Horikawa surface $H(n)$ is a double cover of
 $\mF_{n-3}$ branched over a smooth curve $D_{n-2}$ in the linear
 system $|6C_0+(4n-8)f|$.
 R.~Fintushel and R.~Stern showed that Horikawa surface can be
 decomposed into
 \[H(n)=B_{n-2}\cup D_{n-2} \cup B_{n-2}, \]
 where $B_{n-2}$ is the complement of a neighborhood of the pair
 of $2$-spheres $(C_0+(n-2)f)$ and $C_0$ in $\mF_{n-3}$ (Lemma 2.1
 in~\cite{FS}), and they proved that an elliptic surface $E(n)$ is
 obtained from $H(n)$ by replacing two rational balls $B_{n-2}$ with
 two configurations $C_{n-2}$ (Lemma 7.3 in~\cite{FS}).
 In other words, R.~Fintushel and R.~Stern proved that Horikawa
 surface $H(n)$ can be obtained from an elliptic surface $E(n)$
 by rationally blowing-down two disjoint configurations $C_{n-2}$
 lying in $E(n)$ in smooth category.
 The aim of this article is to prove that the rational blow-down surgery
 above can be performed in complex category, which is following.

\begin{proof}[Proof of Theorem~\ref{thm-main}]
 Note that $\tX_n$ is a double covering of $\mF_n$ branched over $D$,
 and the minimal resolution of two rational double points of type
 $A_1$ and $A_{2n-9}$ in $\tX_n$ is $E(n)$, which is also a double
 cover of $Z_n$ branched over the proper transform of $D$.
 Since the proper transform of $D$ does not meet the contracted linear
 chain of $\mP^1$'s, we have a double cover of $Y_n$ branched
 over the image of the proper transform of $D$ by the map $\psi$.
 We denote this surface by $X_n$.
 Then $X_n$ is the singular surface obtained by contracting two disjoint
 configurations $C_{n-2}$ from an elliptic surface $E(n)$ and
 it has two quotient singularities of class $T$,
 both are of type $\frac{1}{(n-2)^2}(1, n-3)$.
 By the fact that the proper transform of $D$ is disjoint from
 the contracted liner chain of $\mP^1$'s and Lemma~\ref{lem-3.4},
 the map from $X_n$ to $Y_n$ is flat and $H^1(Y_n,\cO_{Y_n}(\bar
 D_{Y_n}))=0$, where $\bar D_{Y_n}$ is the image of $\Delta$ in $Y_n$
 under the contraction $C_{n-2}$ of the liner chain of $\mP^1$'s.
 Therefore we have the following commutative diagram of maps

\[\begin{array}{ccccc}
   \tX_n &  \leftarrow & E(n) & \rightarrow & X_n \\
   \downarrow & & \downarrow & & \downarrow\\
   \mF_n &  \stackrel{\pi}{\leftarrow} & Z_n & \stackrel{\psi}{\rightarrow}
    & Y_n
 \end{array} \]

\n where all vertical maps are double coverings.
 Then, by Theorem~\ref{thm-2.2} and Proposition~\ref{pro-3.2} above,
 the singular surface $X_n$ has a $\mQ$-Gorenstein
 smoothing of two quotient singularities all together.

 Finally, by applying the standard arguments about Milnor fibers
 (\S5 in~\cite{LW} or \S1 in~\cite{Man01}),
 we know that a general fiber of a $\mQ$-Gorenstein smoothing of $X_n$
 is diffeomorphic to the $4$-manifold obtained by rational blow-down
 of $E(n)$.
 And we also know that $H(n)$ has one deformation class
 (\cite{BHPV}, Chapter VII).
 Therefore a general fiber of a $\mQ$-Gorenstein smoothing of $X_n$
 is a Horikawa surface $H(n)$ in complex category.
\end{proof}

\b
\b

\end{document}